\documentclass{ifacconf}

\usepackage{graphicx}      
\usepackage{natbib}        
\usepackage{algorithm,algorithmicx, algpseudocode}
\usepackage{amsmath}
\usepackage{amssymb}
\usepackage{color}
\usepackage{changes}
\newtheorem{assumption}{Assumption}
\newtheorem{definition}{Definition}
\newtheorem{proposition}{Proposition}
\newtheorem{remark}{Remark}
\newtheorem{lemma}{Lemma}

\DeclareMathOperator*{\argmin}{argmin}  
\DeclareMathOperator*{\argmax}{argmax}  
\begin{document}
\begin{frontmatter}

\title{Rate analysis of dual averaging for nonconvex distributed optimization}

\thanks[footnoteinfo]{This work was supported by the Knut and Alice Wallenberg Foundation, the Swedish Foundation for Strategic Research, and the Swedish Research Council.}

\author[First]{Changxin Liu} 
\author[First]{Xuyang Wu} 
\author[First]{Xinlei Yi} 
\author[Second]{Yang Shi} 
\author[First]{Karl H. Johansson}

\address[First]{School of Electrical Engineering and Computer Science, KTH Royal Institute of Technology, and Digital Futures,
	100 44 Stockholm, Sweden (e-mail: \{changxin; xuyangw; xinleiy; kallej\}@kth.se).}

 \address[Second]{Department of Mechanical Engineering, University of Victoria, 
  Victoria, B.C. V8W 3P6, Canada (e-mail: yshi@uvic.ca).}

\begin{abstract}                
This work studies nonconvex distributed constrained optimization over stochastic communication networks.
We revisit the distributed dual averaging algorithm, which is known to converge for convex problems. We start from the centralized case, for which the change of two consecutive updates is taken as the suboptimality measure. 
We validate the use of such a measure by 
showing that it is closely related to stationarity.
This equips us with a handle to study the convergence of dual averaging in nonconvex optimization. We prove that the squared norm of this suboptimality measure converges at rate $\mathcal{O}(1/t)$. Then, for the distributed setup we show convergence to the stationary point at rate $\mathcal{O}(1/t)$. Finally, a numerical example is given to illustrate our theoretical results.

\end{abstract}


\end{frontmatter}

\section{Introduction}

In recent years, distributed optimization has received surged research interests from both academia and industry, because of its capability of delivering high-quality solutions to a system-wide task under the support of a cluster of computing units/agents and real-time communication networks. For a recent overview of distributed optimization, the interested readers are referred to \citep{yang2019survey}.

This work is concerned with the distributed optimization problem where the cost function is the sum of multiple smooth and possibly nonconvex objective functions locally with the agents, the constraint set is common across the agents, and the communication network is time-varying and random.
Such formulation finds wide applications including platooning control of multiple vehicles \citep{shen2022nonconvex}, machine learning \citep{lian2017can}, to name a few. Particularly, the stochastic time-varying communication network is of practical significance because real communication networks suffer from congestion, failure, and random package dropouts. 

Existing works on distributed nonconvex optimization mostly dealt with fixed communication networks; see, e.g., \citep{di2015distributed,hong2017prox,yi2021linear}. 
Recently, \cite{scutari2019distributed,xin2021stochastic,jiang2022distributed} considered nonconvex composite optimization with deterministic time-varying networks. 
However, the communication network is essentially assumed to be connected in every finite steps.
Different from the aforementioned gradient descent based distributed optimization algorithms, distributed dual averaging (DDA) originally proposed by \cite{duchi2011dual} has demonstrated its advantages in simultaneously handling constraints and stochastic communication networks. 
Nevertheless, this type of algorithms were only known to converge for convex problems to the best of our knowledge.

 It is {worth mentioning} that there are a few recent attempts in the literature regarding the convergence of centralized dual averaging (CDA) for nonconvex optimization. For example, \cite{defazio2022adaptivity} established the relation between hyperparameters in CDA and stochastic gradient descent (SGD), and then generalized the analysis in SGD to CDA. However, the analysis is limited to unconstrained optimization, in which SGD and CDA only differ in {the} choice of hyperparameters. However, they may generate distinct trajectories of variables in the presence of constraints \citep{fang2022online}. \cite{heliou2020online} investigated the behavior of dual averaging in online nonconvex optimization with constraints. The authors considered nonsmooth time-varying loss functions with bounded subgradients, which is not applicable to the setup considered in this work.

In this work, we extend the dual averaging based distributed optimization algorithm developed in \citep{liu2022accelerated} to nonconvex constrained problems. The main contributions of this work are as follows. First, we prove the convergence rate of CDA for nonconvex smooth optimization with constraints for the first time. A new measure of suboptimality is defined and its relation to stationarity is discussed. Based on them, we prove the $\mathcal{O}(1/t)$ convergence rate of dual averaging in terms of the squared norm of the suboptimality measure. Then, the results are extended to the distributed setup with stochastic communication networks. Under rather mild conditions, the convergence rate of DDA is proved to be $\mathcal{O}(1/t)$.

\textit{Notation} Given a convex set $\mathcal{X}\subset \mathbb{R}^m$, we denote the normal cone to $\mathcal{X}$ at $x$ by $\mathcal{N}_{\mathcal{X}}(x)=\{ g\in\mathbb{R}^m:  \langle g, y-x \rangle \leq 0, \forall y\in\mathcal{X} \}$. For a real-valued random vector $x$, we define $\lVert x \rVert_{\mathbb{E}} = \sqrt{\mathbb{E}\lVert x \rVert^2 }$. We use $\rho(\cdot)$ to denote the spectral radius of a matrix. A differentiable function $d$ is said strongly convex with modulus $a^{-1}>0$ if 
\begin{equation*}
    d(y) \geq d(x)  + \langle \nabla d(x) , y-x \rangle + \frac{1}{2a} \lVert y -x \rVert^2,  \,\, \forall x, y.
\end{equation*}


\section{Problem Statement}

\subsection{Optimization problem}
Consider the finite-sum constrained optimization problem
\begin{equation} \label{eq:OPT}
	\min_{x\in\mathcal{X}} \left\{ f(x) :=\frac{1}{n} \sum_{i=1}^n  f_i(x)  \right\}
\end{equation}
where each $f_i$ is a smooth and possibly nonconvex function, and $\mathcal{X} \subset \mathbb{R}^m$ is a compact convex set. 
The optimal objective value is denoted as $f^*>-\infty$.  

\begin{assumption}
	\label{assum:smoothness}
	Each $f_i$ is continuously differentiable on an open set that contains $\mathcal{X}$, and $\nabla f_i$ is Lipschitz continuous with Lipschitz constant $L>0$, i.e., 
	\begin{equation*}
		\lVert \nabla f_i(x) - \nabla f_i(y) \rVert \leq L\lVert  x - y \rVert, \quad \forall  x,y\in \mathcal{X}.
	\end{equation*}

\end{assumption}
A direct consequence of Assumption \ref{assum:smoothness} is
\begin{equation*}
	f_i(y) \leq f_i(x) + \langle \nabla f_i(x), y-x \rangle + \frac{L}{2}\lVert y-x \rVert^2 , \quad \forall  x,y\in \mathcal{X}.
\end{equation*}

 For Problem \eqref{eq:OPT}, we recall the stationarity condition, which is a necessary local optimality condition \cite[Definition 3.73]{beck2017first}.
\begin{definition}[stationary point]
   A point $x^*\in\mathcal{X}$ is a stationary point of Problem \eqref{eq:OPT} if $
    -\nabla f(x^*) \in \mathcal{N}_{\mathcal{X}}(x^*).
$
\end{definition}

\subsection{Communication network}
Consider the standard distributed optimization setup, where each agent $i$ holds a local objective function $f_i$ and is only able to communicate with other agents if they are connected in the communication network. At time $t$, a doubly stochastic matrix $P^{(t)}$ is used to describe the network topology and the weights of connected links. In this work, we consider a general setting of stochastic communication networks, i.e., $P^{(t)}$ is a random matrix for every $t$. We denote by $p_{ij}^{(t)}$ the $(i,j)$-th element in $P^{(t)}$. $p_{ij}^{(t)}>0$ only if the two agents $i$ and $j$ are neighbors at $t$. The set of $i$'s neighbors at time $t$ is denoted as $\mathcal{N}_i^{(t)}$.



\begin{assumption}
	\label{graphconnected}
	For every $t\geq 0$, it holds that
		i) $P^{(t)}\mathbf{1} = \mathbf{1}$ and $\mathbf{1}^TP^{(t)} = \mathbf{1}^T$, where $\mathbf{1}$ denotes the all-one vector of dimensionality $n$;
		ii) $P^{(t)}$ is independent of the random events that occur up to time $t-1$; and
		iii) there exists a constant $\beta\in(0,1)$ such that
		\begin{equation}
			\label{eq:sigma-bound}
			\sqrt{\rho\left(\mathbb{E}_t\left[{P^{(t)}}^TP^{(t)}\right] - \frac{\mathbf{1}\mathbf{1}^T}{n}\right)} \leq \beta,
		\end{equation}
		where
  the expectation $\mathbb{E}_t[\cdot]$ is taken with respect to the distribution of $P^{(t)}$ at time $t$.
	\end{assumption}
	
Assumption \ref{graphconnected} is satisfied by a host of common stochastic networks, e.g., randomized gossip \citep{boyd2006randomized} and Bernoulli stochastic networks \citep{kar2008sensor}. Different from the deterministic time-varying networks considered in \citep{nedic2017achieving,xin2021stochastic,jiang2022distributed}, Assumption \ref{graphconnected} does not require the communication network to be connected every finite time steps. In fact, for stochastic {networks defined in} Assumption \ref{graphconnected}, it is possible that the mixing matrix $P^{(t)}$ never produces a deterministic contraction property in finite steps. Thus, the convergence analysis for deterministic time-varying networks cannot be applied or easily extended to the setting of stochastic networks.

This work focuses on the theoretical convergence properties of dual averaging algorithms for nonconvex optimization in both centralized and distributed settings.

\section{Dual Averaging Algorithm for Nonconvex Optimization}

In this section, we present the CDA algorithm and derive its convergence rates for nonconvex problems.

Given constant $a>0$
and an arbitrary variable $x^{(0)} \in \mathcal{X}$, we define a class of proximal functions $d:\mathbb{R}^m\rightarrow \mathbb{R}$.
\begin{definition}[proximal function]
	$d$ is called a proximal function
	if: i) $x^{(0)}$ is the ``$d(\cdot)$-center" of $\mathcal{X}$, i.e.,
$
		x^{(0)} = \argmin_{x\in\mathcal{X}} d(x)$ and $ d(x^{(0)}) = 0
$; ii) $d(x)$ is $a^{-1}$-strongly convex and differentiable.
\end{definition}

Associated with $d$, we define the convex conjugate
\begin{equation*}
	d^*(z) = \max_{x\in\mathcal{X}}\{ \langle z, x \rangle - d(x)\}.
\end{equation*}
According to Danskin's Theorem \citep[Proposition 6.1.1]{bertsekas1999nonlinear}, it holds that
\begin{equation*}
	\nabla d^*(z) = \argmax_{x\in\mathcal{X}}\{ \langle z, x \rangle - d(x)\}.
\end{equation*}

Starting from $x^{(0)}$, CDA produces a sequence of variables $\{ x^{(t)}\}_{t\geq 0}$ according to
\begin{equation}\label{eq:cda_iterate}
	x^{(t)} = \nabla d^* (- {z^{(t)}})
\end{equation}
where 
\begin{equation}\label{eq:cda_iterate_z}
	z^{(t)} = {\sum_{\tau=0}^{t-1} \nabla f(x^{(\tau)})}.
\end{equation}

To investigate the convergence of CDA for nonconvex optimization, we define the following mapping that can be taken as a generalization of the notion of the gradient. The convergence of the proximal gradient descent algorithm \cite[Definition 10.5]{beck2017first} relies on a similar concept, in the sense that they both represent the change of two consecutive updates. Nevertheless, CDA and the proximal gradient descent generally lead to different trajectories for constrained problems. Therefore, the properties of gradient mapping in CDA need to be re-examined.

\begin{definition}[gradient mapping]\label{def:grad_map}
	Suppose that Assumption \ref{assum:smoothness} holds. For any primal-dual pair $({x}^{(t)}, {z}^{(t)})$ generated by \eqref{eq:cda_iterate} and \eqref{eq:cda_iterate_z}, the gradient mapping is defined by
	\begin{equation}\label{eq:grad_map}
		\begin{split}
		G_{a}({x}^{(t)}, { z}^{(t)}) =  \frac{1}{a} \left( \nabla d^*(-{ z}^{(t)}) - \nabla d^*(-z^{(t)}-\nabla{f}({ x}^{(t)})) \right).
	\end{split}
	\end{equation}
\end{definition}
When $\mathcal{X} = \mathbb{R}^m$ and $d(x) = \lVert  x-x^{(0)}\rVert^2/(2a)$, $G_a(x^{(t)},z^{(t)}) = \nabla f(x^{(t)})$ for all $t\geq 0$. In this case, clearly, ${x}^*$ is a stationary point of Problem \eqref{eq:OPT} if and only if there exists $z^*$ such that $G_a(x^*,z^*)=0$. 
In {the unconstrained case}, the relation between $x^{(t)}$ and $z^{(t)}$ is bijective. However, in the presence of constraint, it only holds that \citep[Theorem~26.5]{rockafellar1970convex}:
\begin{equation*}
-z^{(t)}\in \{\nabla d(x^{(t)})\} \oplus \mathcal{N}_{\mathcal{X}}(x^{(t)}), \,\, x^{(t)} = \nabla d^*(-z^{(t)}) \quad  \forall t\geq 0,
\end{equation*}
where $\oplus$ denotes the Minkowski sum defined by
\begin{equation*}
    \mathcal{A} \oplus \mathcal{B}:=\{a+b|a\in \mathcal{A}, b\in\mathcal{B}  \}.
\end{equation*}
\begin{proposition}\label{prop-stationarity}
	${x}^*$ is a stationary point of Problem \eqref{eq:OPT} if and only if there exists some primal-dual pair $(x^*,z^{(t^*)})$ at {some $t^*\ge 0$}, i.e., $x^*= \nabla d^*(-z^{(t^*)})$, such that $G_a(x^{*},z^{(t)})=0, \forall t\geq t^*$.
\end{proposition}

\begin{pf}
\textit{Necessity}. Suppose $x^*$ is a stationary point. 
Pick any $z{'}$ satisfying $x^*= \nabla d^*(-z{'})$, and label the time instant as $t^*$, i.e., $z'=z^{(t^*)}$. By optimality, it holds that
	\begin{equation*}
		\begin{split}
			-\sum_{\tau = 0}^{t^*-1}  \nabla f(x^{(\tau)})  -\nabla d(x^*) &\in  \mathcal{N}_{\mathcal{X}}(x^*).
		\end{split}
	\end{equation*}
 If $x^*$ is a stationary point, we have $ -\nabla f(x^*) \in \mathcal{N}_{\mathcal{X}}(x^*)$ and therefore 
  	\begin{equation*}
		\begin{split}
			-\sum_{\tau = 0}^{t^*}   \nabla f(x^{(\tau)})  -\nabla d(x^*) &\in  \mathcal{N}_{\mathcal{X}}(x^*)  .
		\end{split}
	\end{equation*}
This {together with the strong convexity of $d$} gives us $x^{(t^*+1)} = x^*$ and $G_a(x^{*}, z^{(t^*)})=0$.
By induction, the equality holds for all $t\geq t^*$.

\textit{Sufficiency}. Suppose there exists some $t^*$ such that $G_a(x^{(t)}, z^{(t)})=0, \forall t\geq t^*$. Thus $ x^{(t)}= 
 x^{(t+\tau)}, \forall \tau \geq 1$. Denoting $x^{(t)} = x^*, \forall t\geq t^*$, it holds that
 	\begin{equation*}
		\begin{split}
			v - (t-t^*)  \nabla f(x^{*})   \in  \mathcal{N}_{\mathcal{X}}(x^*) \,\,\mathrm{and}\,\,
			v &\in  \mathcal{N}_{\mathcal{X}}(x^*)   ,\forall t\geq t^*+1
		\end{split}
	\end{equation*}
	where $v = 	-\sum_{\tau = 0}^{t^*-1}  \nabla f(x^{(\tau)})  -\nabla d(x^*)$.
For the sake of contradiction, suppose 
	\begin{equation*}
		-\nabla f(x^*)  \notin  \mathcal{N}_{\mathcal{X}}(x^*).
	\end{equation*}
 Then, there must exist some sufficiently large $t$ such that
 \begin{equation*}
     v - (t-t^*)  \nabla f(x^{*})   \notin  \mathcal{N}_{\mathcal{X}}(x^*),
 \end{equation*}
 which yields a contradiction.
 
\end{pf}

Next, we present the convergence rate of CDA for general nonconvex optimization problems.

\begin{thm}\label{thm:cda}
	Suppose Assumption \ref{assum:smoothness} is satisfied and let $\{x^{(t)}\}_{t\geq 0}$ be the sequence generated by the dual averaging algorithm in \eqref{eq:cda_iterate} and \eqref{eq:cda_iterate_z} with  $a< 2L^{-1}$. Then
	\begin{itemize}
		\item [i)] the sequence $\{f(x^{(t)})\}_{t\geq 0}$ is non-increasing, and
$f(x^{(t)})> \lim_{\tau\rightarrow \infty} f(x^{(\tau)})$ if and only if $x^{(t)}$ is not a stationary point;
		\item [ii)] 	$G_{a}({x}^{(t)}, { z}^{(t)})\rightarrow 0$ as $t\rightarrow \infty$;
		\item [iii)] {for all $k\ge 1$,}
	\begin{equation}\label{eq:convergence_cda}
\min_{t\leq k} \, \lVert G_{a}( x^{(t)}, z^{(t)} )  \rVert^2 \leq \frac{2\left(f(x^{(0)}) - f^*\right)}{  a(2-aL) k}.
	\end{equation}
	\end{itemize}

	
\end{thm}

\begin{remark}
    Theorem \ref{thm:cda} provides a sufficient condition for the parameter $a$, under which CDA converges. In particular, the objective value monotonically decreases until a stationary point is reached. Furthermore, the norm of the suboptimality measure in Definition \ref{def:grad_map} converges to $0$. Finally, the minimum of squared norm of the measure before arbitrary time $k\geq 1$ is bounded from above by $\mathcal{O}(1/k)$.
\end{remark}

\subsubsection{Proof of Theorem \ref{thm:cda}:} Before proving Theorem \ref{thm:cda}, we present Lemma \ref{lem:descent_cda}.
\begin{lemma}\label{lem:descent_cda}
Suppose Assumption \ref{assum:smoothness} holds. For the sequence $\{ x^{(t)}\}_{t\geq 0}$ generated by the dual averaging method in \eqref{eq:cda_iterate} and \eqref{eq:cda_iterate_z}, it holds that
	\begin{equation}\label{eq:cda_descent}
		\begin{split}
		\langle \nabla f(x^{(t)}), x^{(t+1)} - x^{(t)} \rangle 
		\leq  - \frac{1}{a}\lVert  x^{(t+1)} -x^{(t)}\Vert^2.
		\end{split}
	\end{equation}
\end{lemma}

\begin{pf}
    	Let 
	\begin{equation}\label{eq:def_r}
		r_{t+1}(x) := \sum_{\tau=0}^{t} \langle \nabla f(x^{(\tau)}), x - x^{(\tau)}\rangle + d(x),~\forall t\ge 0.
	\end{equation}
    By definition, we have
	\begin{equation}\label{eq:def_r_iter}
		r_{t+1}(x) = r_{t}(x)  + \langle \nabla f(x^{(t)}), x-x^{(t)} \rangle.
	\end{equation}
	According to \eqref{eq:cda_iterate}, we know that $x^{(t)} = \argmin_{x\in\mathcal{X}} r_{t}(x)$. Thus, upon using $r_{t}$ is strongly convex with modulus $1/a$, we obtain
	\begin{equation*}
		r_{t}(x) - r_{t}(x^{(t)}) \geq \frac{1}{2a}\lVert  x -x^{(t)}\Vert^2, \quad \forall x\in \mathcal{X}.
	\end{equation*}
	Taking $x= x^{(t+1)}$, we have
	\begin{equation*}
		\begin{split}
			0&\leq r_{t}(x^{(t+1)}) - r_{t}(x^{(t)}) -\frac{1}{2a}\lVert  x^{(t+1)} -x^{(t)}\Vert^2\\
			& =  r_{t+1}(x^{(t+1)}) - \langle \nabla f(x^{(t)}), x^{(t+1)}-x^{(t)}\rangle  \\
			& \quad - r_{t}(x^{(t)}) -\frac{1}{2a}\lVert  x^{(t+1)} -x^{(t)}\Vert^2,
		\end{split}
	\end{equation*}
	where the equality is due to \eqref{eq:def_r_iter}.
	Therefore,
	\begin{equation}\label{eq:cda_descent_part1}
		\begin{split}
			&	\langle \nabla f(x^{(t)}), x^{(t+1)}-x^{(t)}\rangle  \\
			&\leq   r_{t+1}(x^{(t+1)}) - r_{t}(x^{(t)})-\frac{1}{2a}\lVert  x^{(t+1)} -x^{(t)}\Vert^2. 
		\end{split}
	\end{equation}
	Finally, we use $x^{(t+1)} = \argmin_x r_{t+1}(x)$ and \eqref{eq:def_r} to get
	\begin{equation}\label{eq:rt_lower_bound}
		\begin{split}
			r_{t+1}(x^{(t+1)}) + \frac{1}{2a} \lVert x^{(t+1)} - x^{(t)} \rVert^2 
			\leq r_{t+1}(x^{(t)}) = r_{t}(x^{(t)})   ,
		\end{split}
	\end{equation}
	which together with \eqref{eq:cda_descent_part1} yields \eqref{eq:cda_descent}.
\end{pf}

We are now in a position to prove Theorem  \ref{thm:cda}.

i)	By Assumption \ref{assum:smoothness}, we have
	\begin{equation}\label{eq:Lipschitz_gradient}
		\begin{split}
		&	f(x^{(t+1)})-f(x^{(t)}) \\
			& \leq   \langle \nabla f(x^{(t)}), x^{(t+1)}- x^{(t)} \rangle +  \frac{L}{2}\lVert x^{(t+1)} - x^{(t)} \rVert^2  .
		\end{split}
	\end{equation}
	Using Lemma \ref{lem:descent_cda}, we obtain
	\begin{equation}\label{eq:descent_da}
		\begin{split}
		f(x^{(t)}) -f(x^{(t+1)}) \geq&  \left( \frac{1}{a}-\frac{L}{2}\right) \lVert x^{(t+1)}- x^{(t)} \rVert^2 \\
		=&  \frac{a(2-aL)}{2} \lVert G_{a}( x^{(t)}, z^{(t)} )  \rVert^2,
	\end{split}
	\end{equation}
which implies $f(x^{(t)})  \geq f(x^{(t+1)})$. 
Because the sequence $\{f(x^{(t)})\}_{t\geq 0}$ is non-increasing and bounded from below, it converges.
If $x^{(t)}$ is not a stationary point, then $\sum_{\tau=t}^{\infty}\lVert G_{a}( x^{(\tau)}, z^{(\tau)} )  \rVert^2 \neq 0$ according to Proposition \ref{prop-stationarity}, and therefore $f(x^{(t)})> \lim_{\tau\rightarrow \infty} f(x^{(\tau)})$. If $x^{(t)}$ is a stationary point, then $\sum_{\tau=t}^{\infty}\lVert G_{a}( x^{(\tau)}, z^{(\tau)} )  \rVert^2 {{=}} 0$ and $x^{(\tau)}=x^{(t)}, \forall \tau \geq t$, and thus $f(x^{(t)}) = \lim_{\tau\rightarrow \infty} f(x^{(\tau)})$.

ii) Because the sequence $\{f(x^{(t)})\}_{t\geq 0}$ converges, $f(x^{(t)})- f(x^{(t+1)})$ converges to $0$ as $t\rightarrow \infty$, which in conjunction with \eqref{eq:descent_da} gives the desired result.

iii)
	Summing \eqref{eq:descent_da} over $t=0,1,\dots,k$ yields
	\begin{equation*}
		\begin{split}
	&	\frac{a(2-aL)}{2} \sum_{t=0}^k \lVert G_{a}( x^{(t)}, z^{(t)} )  \rVert^2 \leq f(x^{(0)}) - f(x^{(k+1)}) \\
		&\leq f(x^{(0)}) - f^*
	\end{split}
	\end{equation*}
	where the last inequality follows from $f(x^{(k+1)}) \geq f^*$. Thus \eqref{eq:convergence_cda} holds.

\section{Distributed Dual Averaging for Nonconvex Optimization}

In this section, we revisit the DDA algorithm in \citep{duchi2011dual,liu2022accelerated} (particularly \citep{liu2022accelerated} where smooth problems are considered), and provide its rate of convergence for nonconvex optimization problems in the form of \eqref{eq:OPT}.

{The design of DDA is motivated in \citep{liu2022accelerated}, where the idea is to use dynamic averaging consensus to estimate $z^{(t)}$ in \eqref{eq:cda_iterate_z} in a distributed manner, followed by a similar step to \eqref{eq:cda_iterate} locally performed by each agent with an inexact version of $z^{(t)}$.
The DDA algorithm is detailed in Algorithm \ref{DP-DDA}. First, each agent initializes the algorithm by setting the local variables $x_i^{(0)}$, $z_i^{(0)}$, and $s_i^{(0)}$ properly. At each time $t\geq 1$, each agent exchanges the variables $z_i^{(t-1)}$, $s_i^{(t-1)}$ with its neighbors at time $t-1$, and then computes $z_i^{(t)}$, $x_i^{(t)}$, and $s_i^{(t)}$ according to steps 3--6.}
	
	
	\begin{algorithm}[tb]
		\caption{DDA}
		\label{DP-DDA}
		\begin{algorithmic}[1]
			\Statex {\bfseries Input:} $a>0$, a continuously differentiable and $a^{-1}$-strongly convex proximal function $d$, $x^{(0)}$
			\Statex {\bfseries Output:} {${x}_i^{(t)}, t= 1,2,\dots$}
			\State {\bfseries Initialize:} set $x_i^{(0)}=x^{(0)}$, $z_i^{(0)} = 0$, and $s_i^{(0)} = \nabla f_i(x^{(0)})$ for all $i =1,\dots,n$ 
			\For{$t=1,2,\dots$, each agent $i$ synchronously}
			\State collect $z_j^{(t-1)}$ and $s_j^{(t-1)}$ from all agents $j\in \mathcal{N}_i^{(t-1)}$
			\State update $z_i^{(t)}$ by 
			\begin{equation*}
				z_i^{(t)} = \sum_{j \in \mathcal{N}_i^{(t-1)} \cup \{i\}} p_{ij}^{(t-1)}\left(z_j^{(t-1)} + s_j^{(t-1)} \right)
			\end{equation*}
			\State compute $x_i^{(t)}$ by 
			\begin{equation*}
				x_i^{(t)} = \nabla d^*_t(-z_i^{(t)})
			\end{equation*}
			\State update $s_i^{(t)}$ by 
			\begin{equation*}
				s_i^{(t)} = \sum_{j \in \mathcal{N}_i^{(t-1)} \cup \{i\}} p_{ij}^{(t-1)} s_j^{(t-1)} + \nabla f_i(x_i^{(t)}) - \nabla f_i(x_i^{(t-1)})
			\end{equation*}
			\EndFor
		\end{algorithmic}
	\end{algorithm}



	\subsection{Analysis setup }
	
	Similar to \citep{duchi2011dual,liu2022accelerated}, we construct a sequence of auxiliary variables $\{y^{(t)}\}_{t\geq 1}$ by
	\begin{equation}\label{eq:auxiliary_variable}
		y^{(t)} = \nabla d^*(-\overline{z}^{(t)}),
	\end{equation}
	where $\overline{z}^{(t)} = n^{-1}\sum_{i=1}^n z_{i}^{(t)}$, and $y^{(0)} = x^{(0)}$. For each $x_i^{(t)}, i=1,\dots, n$ and $y^{(t)}$, we have the following relation \cite[Lemma 5]{duchi2011dual}.
	\begin{lemma}\label{lem:xy_to_z}
		For every $t\geq 0$ and $i=1,\dots,n$, there holds
		\begin{equation*}
			\lVert  x_i^{(t)} - y^{(t)} \rVert \leq a{\lVert  z_i^{(t)}- \overline{z}^{(t)}  \rVert}.
		\end{equation*}
	\end{lemma}
	
	To proceed, we recall the analysis from \citep{liu2022accelerated} in quantifying $\lVert  z_i^{(t)}- \overline{z}^{(t)}  \rVert$. First, we introduce the {notations}:
	\begin{equation*}
		\bf{M} = \begin{bmatrix}
			\beta & \beta \\
			aL(\beta+1) & \beta(aL+1)
		\end{bmatrix}
	\end{equation*}
	and
	\begin{equation}\label{eq:notation-1}
		{\bf x}^{(t)}=\begin{bmatrix}
			x_{1}^{(t)}\\ \vdots \\ x_{n}^{(t)}
		\end{bmatrix}, \quad {\bf y}^{(t)}=\begin{bmatrix}
		y^{(t)}\\ \vdots \\ y^{(t)}
	\end{bmatrix}, \quad \overline{g}^{(t)} = \frac{1}{n}\sum_{i=1}^n \nabla f_i(x_i^{(t)}).
	\end{equation}
For the dual variable $\overline{z}^{(t)}$ in \eqref{eq:auxiliary_variable}, one can verify from steps 4 and 6 in Algorithm \ref{DP-DDA} that
\begin{equation*}
    \overline{z}^{(t)} = \overline{z}^{(t-1)} +  \overline{s}^{(t-1)} = \overline{z}^{(t-1)} +  \overline{g}^{(t-1)},
\end{equation*}
where $\overline{s}^{(t)}={n}^{-1}\sum_{i=1}^{n}s_i^{(t)}$.
Next, we remark that the update of $\{y^{(t)}\}_{t\geq 1}$ in \eqref{eq:auxiliary_variable} can be viewed as dual averaging with inexact gradients, whose convergence property is summarized in the following lemma.
	 \begin{lemma}\label{lem:conv_y}
	Suppose Assumption \ref{assum:smoothness} holds. For $y^{(t)}, t= 1,\dots, $ generated by \eqref{eq:auxiliary_variable}, it holds that $\forall  \epsilon>0$
	 	\begin{equation} \label{eq:descent}
	 		\begin{split}
	 			n  \left( f(y^{(t)}) -  f(y^{(t-1)}) \right) 
	 			\leq & \left(\frac{L+\epsilon}{2}-\frac{1}{a}\right) \lVert {\bf y}^{(t)} - {\bf y}^{(t-1)} \rVert^2 \\
	 			&+ \frac{L^2}{2\epsilon} \lVert {\bf y}^{(t-1)}- {\bf x}^{(t-1)}  \rVert^2.
	 		\end{split}
	 	\end{equation}
	 \end{lemma}
	 \begin{pf}
	 	By Assumption \ref{assum:smoothness}, it holds that
	 	\begin{equation}\label{eq:Lipschitz_gradient}
	 		\begin{split}
	 			& f(y^{(t)}) -  f(y^{(t-1)}) \\
	 			&\leq \langle \nabla f(y^{(t-1)}), y^{(t)}- y^{(t-1)} \rangle +  \frac{L}{2}\lVert y^{(t)} - y^{(t-1)} \rVert^2 \\
	 			&=  \langle \nabla f(y^{(t-1)})- \overline{g}^{(t-1)}, y^{(t)}-y^{(t-1)}  \rangle \\
	 			& \quad +  \langle  \overline{g}^{(t-1)}, y^{(t)}-y^{(t-1)}   \rangle +  \frac{L}{2}\lVert y^{(t)} - y^{(t-1)} \rVert^2.
	 		\end{split}
	 	\end{equation}
	 	Using Lemma \ref{lem:descent_cda} over the sequence $\{y^{(t)}\}_{t\geq 1}$ generated by \eqref{eq:auxiliary_variable}, we obtain
	 	\begin{equation*}
	 		\begin{split}
	 		\langle \overline{g}^{(t-1)}, y^{(t)} - y^{(t-1)} \rangle  \leq  - \frac{1}{a}\lVert  y^{(t)} -y^{(t-1)}\Vert^2 .
 		\end{split}
	 	\end{equation*}
	 	Thus
	 	\begin{equation}\label{eq:lem_desc}
	 		\begin{split}
	 				f(y^{(t)})-f(y^{(t-1)})
	 			&\leq   \langle \nabla f(y^{(t-1)})- \overline{g}^{(t-1)}, y^{(t)}-y^{(t-1)}  \rangle   \\
	 			& \quad  -\left(\frac{1}{a} -\frac{L}{2}\right) \lVert y^{(t)}- y^{(t-1)} \rVert^2.
	 		\end{split}
	 	\end{equation}
	 	In addition, {we have
	 		\begin{equation*}
	 			\begin{split}
	 				& \lVert  \nabla f(y^{(t-1)})- \overline{g}^{(t-1)} \rVert^2 \\
	 				&= \left\lVert n^{-1} \sum_{i=1}^n \nabla f(y^{(t-1)})- \nabla f_i(x_i^{(t-1)}) \right \rVert^2 \\ &\leq \frac{1}{n} \sum_{i=1}^n\lVert \nabla f(y^{(t-1)})- \nabla f_i(x_i^{(t-1)})  \rVert^2 \\
	 				& \leq \frac{L^2}{n} \lVert {\bf y}^{(t-1)}- {\bf x}^{(t-1)}  \rVert^2,
	 			\end{split}
	 		\end{equation*}
	 		where the first inequality is due to the convexity of norm square and Jensen's inequality, and the second inequality follows from Assumption \ref{assum:smoothness}.
	 		Therefore, it holds that
	 		\begin{equation}\label{eq:lem_inexact}
	 			\begin{split}
	 				&   \langle \nabla f(y^{(t-1)})- \overline{g}^{(t-1)}, y^{(t)}-y^{(t-1)}  \rangle  \\
	 				& \leq      \frac{\epsilon}{2}\lVert y^{(t)} - y^{(t-1)} \rVert^2 + \frac{1}{2\epsilon} \lVert \nabla f(y^{(t-1)})- \overline{g}^{(t-1)} \rVert^2 \\
	 				& \leq   \frac{\epsilon}{2}\lVert y^{(t)}\! -\! y^{(t-1)} \rVert^2 \!+\! \frac{L^2}{2n\epsilon} \lVert  {\bf y}^{(t-1)}\!- {\bf x}^{(t-1)}\rVert^2 ,
	 				\forall  \epsilon>0.
	 			\end{split}
	 		\end{equation}
	 		Combining \eqref{eq:lem_desc} and \eqref{eq:lem_inexact} completes the proof.} 
	 	
	 	
	 \end{pf}

	

	\subsection{Rate of convergence}
	
	We begin by defining the consensual stationary point in the distributed case. 
 
 \begin{definition}[consensual stationary point]
			A vector $\mathbf{x}^{*} = [x_1^{*};\dots; x_n^{*}]$ is called a stationary solution if 
\begin{equation}\label{eq:residual}
x_1^* = \cdots = x_n^*  \,\, \mathrm{and}  \,\, -\nabla f(x_i^*) \in \mathcal{N}_{\mathcal{X}}(x_i^*), \, \forall i =1, \cdots ,n.
			\end{equation}
		\end{definition}


A sufficient condition to \eqref{eq:residual} is that there exists a primal-dual pair $(y^*,\overline{ z}^{(t^*)})$ at time $t^*$, i.e., $y^*=\nabla d^*(-\overline{z}^{(t^*)})$, such that \begin{equation}\label{eq:residual}
			 n \lVert G_a({y}^{*}, \overline{ z}^{(t)}) \rVert^2 +	 \sum_{i=1}^n\left \lVert {x }_i^{*} - y^{*} \right\rVert^2   =0, \,\, \forall t\geq t^*
			\end{equation}
   where $G_a(\cdot,\cdot)$ is defined in \eqref{eq:grad_map}.
To see this, we note that \eqref{eq:residual} implies
			\begin{equation*}
				\lVert G_a({y}^{*}, \overline{ z}^{(t)}) \rVert^2 =0 \quad \mathrm{and} \quad  \left \lVert {x }_i^{*} - y^{*} \right\rVert^2  =0 \,\, \forall i=1,\dots,n
			\end{equation*}
    for all $t \geq t^*$,
			where the former ensures that $y^{*}$ is a stationary point in the centralized case, and the latter gives $x_i^* = y^*, \forall i=1,\dots,n$. 

	\begin{thm}\label{thm:sublinear}
		Suppose Assumptions \ref{assum:smoothness} and \ref{graphconnected} are satisfied.
  If the constant $a$ satisfies 
		\begin{equation*}
			\frac{1}{a} > L \cdot \max \left\{ \frac{1}{2}+\frac{4}{3(1-\rho({\bf M}))} ,  \frac{2\beta }{(1-\beta)^2} \right\},
		\end{equation*}

		then, for all $t\geq 1$, it holds that {
  \begin{equation}\label{thm-eq:residual-inf}
      \lim_{t\rightarrow +\infty} n \lVert G_a({ y}^{(t-1)}, \overline{ z}^{(t-1)}) \rVert^2_{\mathbb{E}} \!+\! \lVert  {\bf x}^{(t-1)} - {\bf y}^{(t-1)}\rVert^2_{\mathbb{E}} = 0
  \end{equation}}
  and
		\begin{equation}\label{thm-eq:residual}
   \min_{\tau \leq t}  n \lVert G_a({ y}^{(\tau-1)}, \overline{ z}^{(\tau-1)}) \rVert^2_{\mathbb{E}} + \lVert  {\bf x}^{(\tau-1)} \!- {\bf y}^{(\tau-1)}\rVert^2_{\mathbb{E}} \leq\! \frac{C}{t}
		\end{equation}
		where 
		\begin{equation*}
		\begin{split}
			C: = &\left( \min \left\{\frac{3L(1-\rho({\bf M}))}{8},  {a} - \frac{a^2L}{2} - \frac{4a^2L}{3(1-\rho({\bf M}))} \right\} \right)^{-1} \\
		&	\times \left(\frac{2 \pi^2}{3 L (1-\rho({\bf M}))} + n\left(  f(y^{(0)}) - f^*\right) \right).
		\end{split}
		\end{equation*}
	\end{thm}
	
	

 \subsubsection{Proof of Theorem \ref{thm:sublinear}:}
	Before proving Theorem \ref{thm:sublinear}, we provide the following lemma. Its proof is similar to the proof of \citep[Lemma 5]{liu2022accelerated}. Due to space limitations, the proof is omitted.
	\begin{lemma}\label{lem:consensus}
		If Assumption \ref{graphconnected} holds and
		\begin{equation}\label{eq:1st_condition}
			{a} <   \frac{(1-\beta)^2} {2\beta L},
		\end{equation}
		then, for ${\bf x}^{(t)}$ and ${\bf y}^{(t)}$ defined in \eqref{eq:notation-1}, it holds that
		\begin{equation*}
			\begin{split}
			\sum_{\tau=1}^t  \lVert {\bf x}^{(\tau)} -{\bf y}^{(\tau)} \rVert^2_\mathbb{E} \leq & \frac{8}{9(1-\rho({\bf M}))^2}\sum_{\tau=0}^{t-1}  \lVert {\bf y}^{(\tau+1)} - {\bf y}^{(\tau)} \rVert^2_\mathbb{E} \\
			&+ \frac{8\pi^2}{9L^2(1-(\rho({\bf M}))^2)} 
		\end{split}
		\end{equation*}
		where 
		$
		\pi^2 =  {\sum_{i=1}^n\left\|\nabla f_i(x^{(0)}) -  \overline{g}^{(0)}\right\|^2}.
		$
	\end{lemma}
	
	Now we are ready to prove Theorem \ref{thm:sublinear}.
	
		Recall \eqref{eq:descent}
		\begin{equation*} 
			\begin{split}
			n  \left( f(y^{(t)}) -  f(y^{(t-1)}) \right) \leq & \left(\frac{L+\epsilon}{2}-\frac{1}{a}\right) \lVert {\bf y}^{(t)} - {\bf y}^{(t-1)} \rVert^2\\
			& + \frac{ L^2}{2\epsilon} \lVert {\bf y}^{(t-1)}- {\bf x}^{(t-1)}  \rVert^2.
		\end{split}
		\end{equation*}
		Summing it from $1$ to $t$ yields
		\begin{equation*}
			\begin{split}
				n\left( f(y^{(t)}) - f(y^{(0)}) \right)
				& \leq    \left(  \frac{L+\epsilon}{2}-\frac{1}{a}  \right) \sum_{\tau=1}^t  \lVert {\bf y}^{(\tau)} - {\bf y}^{(\tau-1)} \rVert^2 \\
				& \quad + \frac{L^2}{2\epsilon} \sum_{\tau=1}^t \lVert {\bf y}^{(\tau-1)}- {\bf x}^{(\tau-1)}  \rVert^2.
			\end{split}
		\end{equation*}
		Taking expectation on both sides, we obtain
		\begin{equation*}
			\begin{split}
				& n \mathbb{E}\left[ f(y^{(t)}) - f(y^{(0)}) \right]  \\
				& \leq \left( \frac{L+\epsilon}{2} -\frac{1}{a} \right) \sum_{\tau=1}^t  \lVert {\bf y}^{(\tau)} - {\bf y}^{(\tau-1)} \rVert_{\mathbb{E}}^2 \\
				& \quad -  \frac{L^2}{2\epsilon} \sum_{\tau=1}^t \lVert {\bf y}^{(\tau-1)}- {\bf x}^{(\tau-1)}  \rVert_{\mathbb{E}}^2+  \frac{8\pi^2}{9\epsilon(1-(\rho({\bf M}))^2)}  \\
				& \quad  +   \frac{8L^2}{9\epsilon(1-\rho({\bf M}))^2}\sum_{\tau=0}^{t-1}  \lVert {\bf y}^{(\tau+1)} - {\bf y}^{(\tau)} \rVert_{\mathbb{E}}^2 \\
				& =  \left(  \frac{L+\epsilon}{2} + \frac{8L^2}{9\epsilon(1-\rho({\bf M}))^2} -\frac{1}{a} \right) \sum_{\tau=1}^t  \lVert {\bf y}^{(\tau)} - {\bf y}^{(\tau-1)} \rVert_{\mathbb{E}}^2 \\
				& \quad - \frac{L^2}{2\epsilon} \sum_{\tau=1}^t \lVert {\bf y}^{(\tau-1)}- {\bf x}^{(\tau-1)}  \rVert^2  + \frac{8  \pi^2}{9\epsilon (1-\rho({\bf M}))^2}.
			\end{split}
		\end{equation*}
		This is equivalent to 
		\begin{equation*}
			\begin{split}
				& a^2\left( \frac{1}{a} - \frac{L+\epsilon}{2} - \frac{8L^2}{9\epsilon(1-\rho({\bf M}))^2} \right) n\sum_{\tau=1}^t  \lVert {\bf y}^{(\tau)} - {\bf y}^{(\tau-1)} \rVert^2_{\mathbb{E}}  \\
				&\quad + \frac{L^2}{2\epsilon} \sum_{\tau=1}^t \lVert  {\bf x}^{(\tau-1)} - {\bf y}^{(\tau-1)}\rVert^2_{\mathbb{E}} \\
				& \leq  \frac{8  \pi^2}{9\epsilon (1-\rho({\bf M}))^2} + n\left(  f(y^{(0)}) - f^*\right) <+\infty
			\end{split}
		\end{equation*}
	because of 
			\begin{equation*}
		\lVert {\bf y}^{(\tau)} - {\bf y}^{(\tau-1)} \rVert^2 = na^2 \lVert G_a({\bf y}^{(\tau-1)}, {\bf z}^{(\tau-1)}) \rVert^2.
	\end{equation*}  
 Thus, \eqref{thm-eq:residual-inf} holds.
		Finally, we set $\epsilon = {4L}/{(3(1-\rho({\bf M})))}>0$ 
	to obtain \eqref{thm-eq:residual}.

\section{Numerical Example}

Consider the distributed principal component analysis (PCA) problem
\begin{equation*}
    \min_{\lVert x \rVert\leq 1} \,\, f(x) :=-\sum_{i=1}^n\lVert M_i x \rVert^2
\end{equation*}
where $n=50$. 
Each agent $i$ possesses a data matrix $M_i\in\mathbb{R}^{30\times 500}$, where each row $M_i^j, j=1,\dots, 30$ is randomly generated with zero mean and $\lVert M_i^{j} \rVert \leq 1$.
For the communication network among agents, we consider the Bernoulli stochastic network \citep{kar2008sensor}, where a complete graph is taken as the supergraph and at each time $t$ every edge of the set of edges of the supergraph is activated with probability $0.1$.
Based on it, a Laplacian-based weight matrix \citep{xiao2005scheme} is used at each time $t$.
    
We initialize Algorithm \ref{DP-DDA} by randomly generating a $500$-dimensional vector with i.i.d. elements drawn from the standard Normal distribution and then projecting it onto the constraint to get $x^{(0)}$. 
Set the parameter $a=1$, and $d(x)= \lVert x-x^{(0)}\rVert^2/2$ accordingly. Comparison is made with the distributed proximal gradient algorithm (DPGA) in \citep{jiang2022distributed}. The stepsize for DPGA is set as $1e{-4}$ in order to stabilize the updates. We remark that DPGA does not have convergence guarantees in stochastic communication networks. 

\begin{figure}
	\centering
	\includegraphics[width=3.4in]{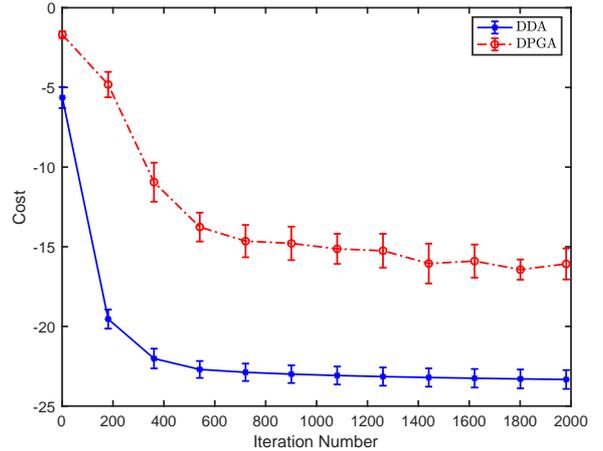}
	\caption{Convergence of the cost $f(x^{(t)})$.}
	
	\label{fig:cost}
\end{figure}

\begin{figure}
	\centering
	\includegraphics[width=3.4in]{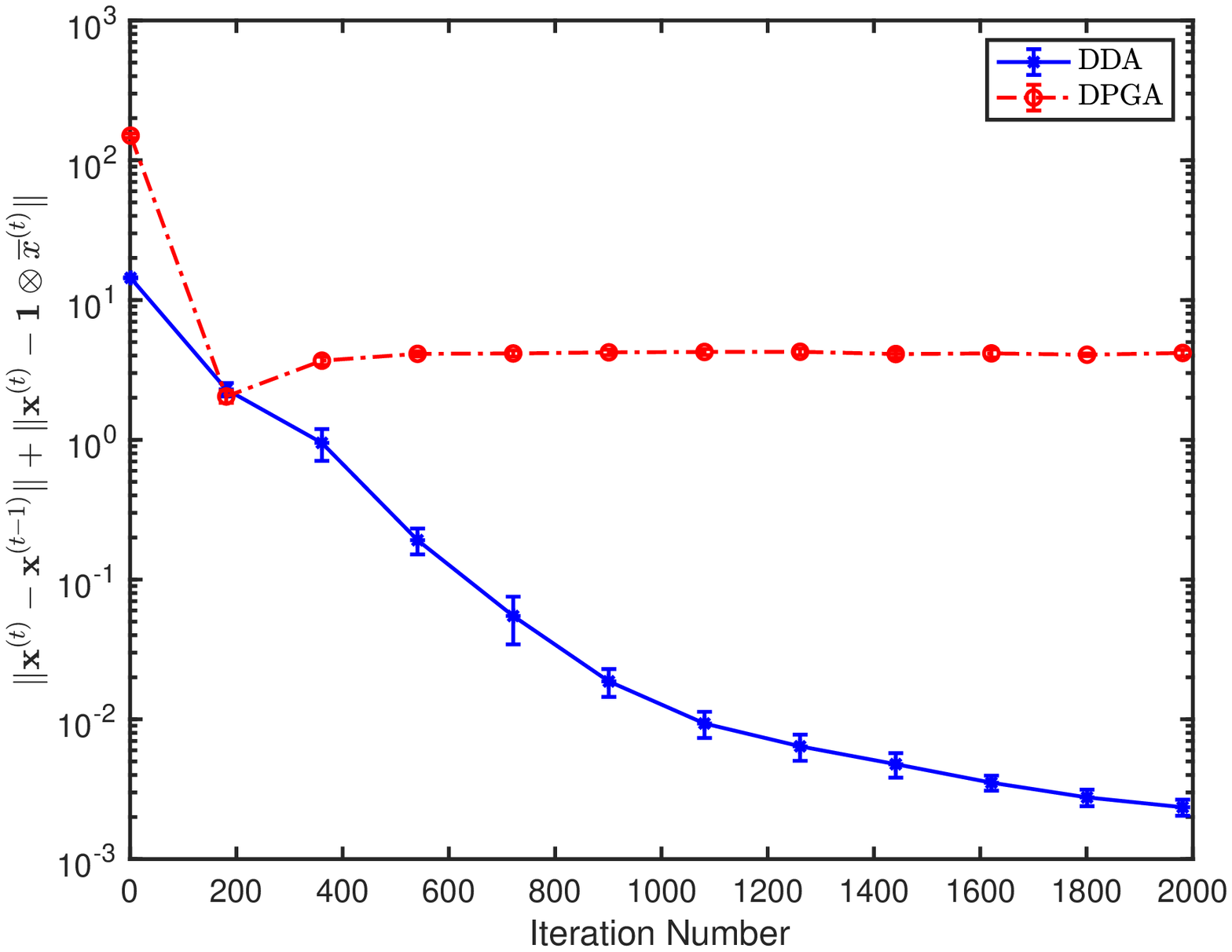}
	\caption{Convergence of $\lVert {\bf x}^{(t)}-{\bf x}^{(t-1)}\rVert + \lVert {\bf x}^{(t)} - {\bf 1} \otimes  \overline{x}^{(t)}  \rVert$. 
	}
	\label{fig:change}
\end{figure}

The experiment was repeated $10$ times with random seeds. We evaluate the performance of the algorithm via the values of the cost function and the sum of difference in two consecutive updates and consensus error, i.e., $\lVert {\bf x}^{(t)}-{\bf x}^{(t-1)}\rVert + \lVert {\bf x}^{(t)} - {\bf 1} \otimes  \overline{x}^{(t)}  \rVert, \,\, t\geq 1$ where $ \overline{x}^{(t)} = n^{-1}\sum_{i=1}^n x_i^{(t)}$ and $\otimes$ denotes the Kronecker product. 
We remark that the latter is an approximation of the residual term in Theorem \ref{thm:sublinear}.
Their mean and standard deviation in $10$ runs by the two algorithms are plotted in Figures \ref{fig:cost} and \ref{fig:change}. Note that a lower value of the cost suggests a closer distance from ${\bf x}^{(t)}$ to the principal eigenvector. From the figures, we observe that for DDA both the cost and the residue converge. In addition, the convergence of DDA is faster than PDGA, since the latter has to use a much smaller stepsize to avoid divergence in this experiment. 
In contrast, DDA remains convergent under a larger range of parameters.
This highlights the advantage of DDA in dealing with stochastic communication networks.

\section{Conclusion}

This work examined the convergence rate of dual averaging for nonconvex constrained smooth optimization problems in both centralized and distributed settings. We developed a new suboptimality measure and established its relation to stationarity. The squared norm of such measure converges at rate $\mathcal{O}(1/t)$ for CDA. Under mild conditions on the stochastic communication network, the rate of DDA is proved to be $\mathcal{O}(1/t)$. For future research, we are interested in speeding up DDA for a special class of nonconvex problems satisfying the Kurdyka-\L{}ojasiewicz condition.


\bibliography{sample}             
                                                   








\end{document}